\documentclass[11pt,reqno]{amsart}
\usepackage[T1]{fontenc}
\usepackage{lmodern}
\usepackage[margin=1in]{geometry}
\usepackage{amsmath,amssymb,amsthm,mathtools}
\usepackage{microtype}
\usepackage{xcolor,booktabs,array,tabularx}
\usepackage[colorlinks=true,linkcolor=blue!40!black,citecolor=blue!40!black,urlcolor=blue!40!black]{hyperref}
\newtheorem{theorem}{Theorem}[section]
\newtheorem{lemma}[theorem]{Lemma}
\newtheorem{proposition}[theorem]{Proposition}
\newtheorem{corollary}[theorem]{Corollary}
\theoremstyle{remark}\newtheorem{remark}[theorem]{Remark}
\theoremstyle{definition}\newtheorem{example}[theorem]{Example}
\theoremstyle{definition}
\usepackage{comment}
\newcommand{\comments}[1]{}
\newcommand{\F}{\mathbb F}
\newcommand{\kk}{\overline{\mathbb F}_p}
\newcommand{\Res}{\operatorname{Res}}
\newcommand{\degt}{\deg_{Y,Z}}
\newcommand{\ord}{\operatorname{ord}}
\numberwithin{equation}{section}
\allowdisplaybreaks[1]
\title[Stepanov's method and sumset estimates]
{Stepanov's method and uniform sumset estimates\\in arbitrary finite fields}
\author{Albert Cochrane}
\address{ Department of Mathematics\\
          Kansas State University\\
          Manhattan, KS 66506}
\email{albertc@ksu.edu}

\author{Todd Cochrane}
\address{ Department of Mathematics\\
          Kansas State University\\
          Manhattan, KS 66506}
\email{cochrane@math.ksu.edu}

\author{Craig Spencer}
\address{ Department of Mathematics\\
          Kansas State University\\
          Manhattan, KS 66506}
\email{cvs@math.ksu.edu}
\subjclass[2020]{Primary 11B30; Secondary 11T30, 11P70}
\keywords{finite field, sum-product, multiplicative subgroup, sumset, Stepanov method}

\begin{document}
%\raggedbottom
%\setlength{\emergencystretch}{1.5em}
\begin{abstract}
Let $A<\mathbb F_q^\times$ be a proper multiplicative subgroup that generates $\mathbb F_q$ over its prime field.
We prove that
\[
 |A+A|\ge \tfrac18|A|^{3/2}
\]
whenever $|A|<q^{2/3}$, uniformly in the characteristic and extension degree, without quantitative restrictions on intersections with proper subfields.
Moreover, in contrast with what is expected in prime fields, we show that no uniform bound with exponent larger than $5/3$ holds in $\mathbb F_q$ even for $|A| \approx q^{1/2}$.
The proof's key idea is to run the prime-field Stepanov method in reverse, constructing an auxiliary polynomial whose specialization vanishes identically instead of requiring it to be nonzero.
\end{abstract}
\maketitle

\section{Introduction}

Let $A$ be a multiplicative subgroup of a finite field.
Since $AA=A$, the size of $A+A$ measures the extent to which a set with complete multiplicative structure must expand under addition.
Subgroups have consequently been interesting examples in finite-field sum-product estimates.

Write
\[
 q=p^f,\qquad A\le \F_q^\times,\qquad t=|A|,
\]
where $p$ is prime.
Over a prime field, the estimate $|A+A|\gg t^{3/2}$ was established by Heath-Brown and Konyagin \cite{HBK} for $t<p^{2/3}$, and stronger exponents are known for smaller $t$.
Over a general field $\mathbb F_q$, the size of $p$ in comparison to $t$ and the presence of proper subfields create additional obstructions.
Little progress has been made for a similar type of lower bound over $\mathbb F_q$; see Table \ref{tab:benchmarks} for a summary of existing estimates over $\mathbb F_p$ and $\mathbb F_q$.
The natural ambient field for this question is the field generated by $A$, so we assume $\F_p(A)=\F_q$.
For $t>1$, this is equivalent to $f=\ord_t(p)$.
This convention excludes the trivial cases where $A$ is contained in a proper subfield, but still allows $A$ to have large intersections with one.
Our main theorem is the following.

\begin{theorem}\label{thm:main}
Let $q=p^f$, and let $A \le \F_q^\times$ be a multiplicative subgroup.
Suppose that
\[
t=|A|,\qquad \F_p(A)=\F_q,\qquad t<(\tfrac 83 q)^{2/3}.
\]
Then
\[
|A+A|\ge \tfrac18\,t^{3/2}.
\]
\end{theorem}

The constant is absolute: there is no restriction on $p$, on $f$, or on the parity of $t$.
In particular, the size condition on $t$ involves $q$, rather than $p$,
and no quantitative subfield-avoidance hypothesis is imposed akin to what one often sees for sum-product estimates in $\mathbb F_q$.

The hypothesis $t\lessapprox q^{2/3}$ is the natural constraint for a three-halves bound since at $t = q^{2/3}$ we have $t^{3/2} = q$.
For $t\gtrapprox q^{2/3}$, classical arguments already give a positive proportion of the field:
\begin{equation}\label{eq:classical-dense}
 |A+A|\ge \frac{q\,t^3}{t^3+q^2}\ge \frac {64}{73}q,
 \quad\text{if}\  t\ge  (\tfrac 83q)^{2/3}.
\end{equation}
The first inequality is a standard consequence of Gauss sums and Parseval's identity.
For completeness, Section~\ref{sec:dense} gives the short argument.
Combining the two ranges gives a convenient all-size formulation.

\begin{corollary}\label{cor:all}
Let $A\le \F_q^\times$ be a multiplicative subgroup with $|A|=t$ and $\F_p(A)=\F_q$.
Then
\[
 |A+A|\ge \min\{\tfrac 18 t^{3/2},\ \tfrac {64}{73} q\}.
\]
\end{corollary}
\begin{proof}
Use Theorem~\ref{thm:main} when $t<(\frac 83 q)^{2/3}$, and \eqref{eq:classical-dense} when $t\ge (\frac 83 q)^{2/3}$.
\end{proof}

Theorem~\ref{thm:main} can also be read as an \textit{inverse theorem}, showing that a subgroup with small doubling must be contained in a subfield of comparable size.
This is similar in spirit to Freiman--Ruzsa-type theorems over finite fields.
To be specific, Ruzsa \cite{Ruzsa1999} with improvement by Even-Zohar and Lovett \cite{EvenZoharLovett2014} established that if $A \subseteq \mathbb{F}_q$,
with $|{}A+A|{}\le K|{}A|{}$, then $A$ is contained in an affine subspace (a translate of an $\mathbb{F}_{p}$-subspace) $V$ of size:
$|{}V|{}\le p^{2K-1}|{}A|{}$. \footnote{for odd $p$ this follows from estimate (5) in the proof of \cite[Theorem 3]{EvenZoharLovett2014}
and for $p = 2$ from \cite[Theorem 12]{ ez2012}. For $p>2$, $K \ge 8$ they prove the stronger bound $|V| \le p^{2K-2}|A|/(2K-1)$.}
Here we obtain:

\begin{corollary}
Let $K$ be a constant, $q$ be a prime power and $A$ be a subgroup of $\mathbb F_q^\times$ (without the assumption $\mathbb F_p(A)=\mathbb F_q)$ such that $|A|>64K^2$ and $|A+A| \le K |A|$.
Then $A$ is contained in a subfield $\mathbb F_Q$ of $\mathbb F_q$ with $Q \le \frac {73}{64}K|A|$.
\end{corollary}
\noindent The subfield $\mathbb F_Q$ is proper if $64K^2<|A|<\frac {64}{73K}\, q$.

\begin{proof}
Suppose that $t>64K^2$ and $|A+A| \le Kt$.
Let $\mathbb F_Q$ be the field generated by $A$, $\mathbb F_Q=\mathbb F_p(A)$.
Then $K < \frac 18 t^{1/2}$ and $|A+A| < \frac 18 t^{3/2}$.
Thus by the preceding corollary, applied with $Q$ in place of $q$, $|A+A| \ge \frac {64}{73}Q$, and so $Q \le \frac {73}{64} Kt$.
\end{proof}

We can ask what the best possible uniform bound on $|A+A|$ might be for small subgroups of $\mathbb F_q^\times$.
Say, for $t\le \sqrt{q}$, can the exponent $3/2$ in Theorem~\ref{thm:main} be improved?
In prime fields, it is possible that we have $|A+A| \gg_\varepsilon |A|^{2-\varepsilon}$ for $|A|\lessapprox p^{1/2}$.
For a general field, the best one can hope for is $|A+A| \gg |A|^{5/3}$ for $|A|\lessapprox q^{1/2}$, as the following example indicates.

\begin{example}
Let $2 \le a<b$ be coprime integers, $H=\mathbb F_{p^a}^\times$, $K=\mathbb F_{p^{b}}^\times$, $q=p^{ab}$, and $A=HK < \mathbb F_q^\times$, so that
\[
|A|=\frac {|H||K|}{|H\cap K|}= \frac {(p^a-1)(p^{b}-1)}{p-1} \approx p^{a+b-1}=q^{\frac 1a+\frac 1b-\frac 1{ab}}.
\]
We claim that $|A+A| \approx |A|^{2-\frac 2{a+b-1}}$.
In particular, for $a=3, b=4$, $A$ is a subgroup of cardinality $|A| \approx \sqrt{q}$ and $|A+A| \approx |A|^{5/3}$.

Let $\{v_1, \dots, v_a\}$, $\{w_1, \dots, w_{b}\}$ be bases for $\mathbb F_{p^a}, \mathbb F_{p^{b}}$ over $\mathbb F_p$, respectively,
and $\{v_iw_j: 1 \le i \le a, 1 \le j \le b\}$ a basis for $\mathbb F_q$.
View the elements $\sum_{i,j} c_{ij}v_iw_j$ of $\mathbb F_q$ as $a$ by $b$ matrices with entries $c_{ij}$, $1 \le i \le a$, $1 \le j \le b$.
A typical element of $A$ is given by $\sum_{i=1}^a x_i v_i \sum_{j=1}^{b} y_jw_j = \sum_{i,j} x_iy_j v_iw_j$ with $\mathbf x \neq \mathbf 0$, $\mathbf y \neq \mathbf 0$.
Thus we can identify $A$ with the set of all rank-1 $a$ by $b$ matrices $[x_iy_j]$ over $\mathbb F_p$.
This is the full set of rank-1 matrices.
The sum of any two such matrices has rank at most 2, and conversely any $a$ by $b$ matrix of rank at most 2 can be expressed as a sum of two rank-1 matrices.
Thus $A+A$ corresponds to the set of all $a$ by $b$ matrices of rank at most 2.
The number of such matrices is roughly obtained by letting the first two row vectors be arbitrary and choosing the remaining $a-2$ row vectors to be arbitrary elements of the linear span of the first two vectors.\footnote{The exact number of  $a$ by $b$ matrices of rank at most 2 is $1+ \frac {(p^a-1)(p^b-1)}{p-1}+\frac {p(p^a-1)(p^{a-1}-1)(p^b-1)(p^{b-1}-1)}{(p^2-1)(p-1)}$.
}
Thus altogether,
\[
|A+A| \approx p^{b}p^{b} p^{2(a-2)} \approx p^{2a+2b-4} \approx |A|^{\frac {2a+2b-4}{a+b-1}}=|A|^{2-\frac 2{a+b-1}}.
\]
\end{example}

\begin{remark}
Regarding the sum-product phenomenon, the preceding example gives a subset $A$ of $\mathbb F_q$ not contained in any proper subfield of $\mathbb F_q$, with $|A| \approx q^{1/2}$ and
\[
\max\{ |A A|,\ | A+A|\} \lessapprox |A|^{5/3}.
\]
For any prime $p$, such a set $A$ exists in $\mathbb F_{p^{12}}^\times$.
For odd primes, by using the quadratic residues in $\mathbb F_{p^a}^\times$ and $\mathbb F_{p^b}^\times$ for $H$ and $K$ instead of the full groups,
one can construct in a similar manner a group $A$ with $|A|$ strictly less than $\sqrt{q}$ having $|A+A|\lessapprox |A|^{5/3}$.

\end{remark}

\begin{remark}
The statement of Theorem \ref{thm:main} and its corollaries hold identically if we replace $A+A$ with a more general sum of cosets.
Thus for any $g \in \mathbb F_q^\times$, we have
\begin{equation}\label{eq:coset}
 |A+gA|\ge
 \min\left\{\tfrac18t^{3/2},\ \tfrac{64}{73}q\right\}
\end{equation}
provided that $\mathbb F_p(A)=\mathbb F_q$.
In particular, this inequality applies to the difference set $|A-A|$.
We outline the modifications needed to prove this generalization in Section \ref{cosetsec}.
\end{remark}

\subsection{Notation}

We use $X\gg Y$ to mean $X\ge cY$ for an absolute positive constant $c$, $X\asymp Y$ when both inequalities hold, and use $L=\log(2t)$ when recording logarithmic factors in such inequalities.

When working with polynomials or rational functions we use upper case $X,Y,Z,...$ to denote indeterminate symbols, and lower case $x,y,z,...$ denote elements of a given field.
Thus $F=F(X,Y,Z)\in \mathbb F_q[X,Y,Z]$ is viewed as a formal object in a polynomial ring, while $F(x,y,z)$ is a function on a given domain.

For a rational function $R=P/Q$ over a given field, with $P,Q$ relatively prime polynomials, we define $\deg R=\max\{\deg P, \deg Q\}$.
Constants have degree zero.

For convenience, if $F=F(X,Y,Z)$ is a polynomial we frequently call $X$ the {\it source variable},
and $Y,Z$ the target variables (since they are often specialized to $X^t,(1-X)^t$) and refer to $\deg_X F$ as the {\it source degree} of $F$ and $\deg_{Y,Z} F$ as the {\it target degree}.
Here, $\deg_{Y,Z} F$ is the total degree of $F$ viewed as a polynomial in $Y$ and $Z$.

\section{Background and outline of proof}

\subsection{Prior estimates and the role of the characteristic}
\label{subsec:literature}

For subgroups of prime fields, Cochrane and Pinner~\cite[Theorem 5.2(b)]{CP} obtained the explicit estimate
\begin{equation}\label{eq:CP}
 |A+A|\ge \tfrac14\,t\bigl(\sqrt{t+1}+1\bigr)
 >\tfrac14\,t^{3/2},
 \qquad t<p^{2/3}.
\end{equation}
using the Stepanov method as developed by Heath-Brown and Konyagin~\cite{HBK}.
%This gives a stronger constant than Theorem~\ref{thm:main}
%for prime fields.

Subsequent work developed stronger prime-field estimates for smaller subgroups, using additive shifts, higher energies,
and incidence methods; see Schoen--Shkredov~\cite{SchoenShkredov}, Vyugin--Shkredov~\cite{VyuginShkredov}, Shkredov~\cite{Shkredov},
Hart~\cite{Hart}, and Murphy--Rudnev--Shkredov--Shteinikov~\cite{MRSS}.

Regarding $\mathbb F_q$, with $q \ge p^2$, there is a three-halves estimate of Roche-Newton, Rudnev, and Shkredov~\cite[Corollary 4]{RNR},
$ |A+A|\gg \min\{t^{3/2},p\}, $ which matches the strength of Corollary \ref{cor:all}, but only for subgroups of very small size, namely $t<p^{2/3}$.
Their argument uses Rudnev's point--plane incidence theorem~\cite{Rudnev}; the characteristic restriction on the incidence populations leads to the threshold $t\ll p^{2/3}$.
There are related estimates for more general sum products $B(B+B)$ with $B$ an arbitrary subset of $\mathbb F_q$.
The collision method of Aksoy Yazici, Murphy, Rudnev, and Shkredov~\cite[Corollary 3(a)]{AMRS} gives $|B(B+B)|\gg |B|^{3/2}$ for $|B|<p^{2/3}$.
For a subgroup, $A(A+A)=A+A$, so this also gives the three-halves growth $|A+A| \gg t^{3/2}$ for $t<p^{2/3}$.
Theorem~\ref{thm:main} replaces the scale $p^{2/3}$ by $q^{2/3}$ for the size of $t$.

\begin{table}[htbp]
\centering
\small
\setlength{\tabcolsep}{4pt}
\renewcommand{\arraystretch}{1.22}
\caption{Selected estimates for $|A+A|$ for $\mathbb F_p$ and $\mathbb F_q$.  Here $L=\log(2t)$,  $t \ge 2$}
\label{tab:benchmarks}
\begin{tabularx}{\textwidth}{@{}>{\raggedright\arraybackslash}p{0.14\textwidth}>{\raggedright\arraybackslash}X >{\raggedright\arraybackslash}p{0.25\textwidth}>{\raggedright\arraybackslash}p{0.13\textwidth}@{}}
\toprule
Ambient field & Hypotheses & Bound for $|A+A|$ & Reference\\
\midrule
$\F_p$& $t<\log_3p$& $\ge \frac 12 t^2+1$ & \cite{CCS2023}\\
$\F_p$& $t\ll p^{1/2}$,\  $t$ even & $\gg t^{5/3}L^{-1/2}$& \cite{VyuginShkredov}\\
$\F_p$& $t\ll p^{5/9}L^{-1/18}$& $\gg t^{8/5}L^{-3/10}$& \cite{Hart}\\
$\F_p$& $t<p^{2/3}$& $>\tfrac14t^{3/2}$& \cite{CP}\\
$\F_q$& Any subgroup &
$\gg\min\{t^{3/2},p\}$& \cite{RNR}\\
%$\F_q$& $\F_p(A)=\F_q,\ t\le q^{3/5}$&
%$\ge t^{7/6}$& \cite{CCPS1}\\
%$\F_q$& $\F_p(A)=\F_q,\ t\le p$&
%$\ge \tfrac 49 t^{4/3}$& \cite{CCPS1}\\
$\F_q$& $\F_p(A)=\F_q,\ t<q^{2/3}$&
$\ge \tfrac 18 t^{3/2}$& Theorem~\ref{thm:main}\\
\bottomrule
\end{tabularx}
\end{table}

General finite-field sum-product theory takes another route by controlling concentration in subfields.
Following the foundational work of Bourgain--Katz--Tao~\cite{BKT},
quantitative estimates in nonprime fields were developed by Katz--Shen~\cite{KatzShen} and Li--Roche-Newton~\cite{LiRN}, among others.
For a precise comparison, Li--Roche-Newton's Theorem 1.4 gives, for $B\subseteq\F_q^\times$,
\[
 \max\{|B+B|,|BB|\}\gg
 |B|^{12/11}(\log(2|B|))^{-5/11}
\]
provided $|B\cap cG|\le |G|^{1/2}$ for every subfield $G\subseteq\F_q$ and every scalar $c$.
Here the hypothesis includes $G=\F_q$, and hence $|B|\le\sqrt q$.
For subgroups, $AA=A$, so this supplies a sumset estimate $|A+A| \gg |A|^{12/11-\varepsilon}$ under the subfield avoidance assumption.

\subsection{Why the prime-field arguments do not transfer directly}
\label{subsec:obstructions}

There are a couple obstructions to simply extending the method of Stepanov \cite{Stepanov} for prime fields,
as developed by Heath-Brown and Konyagin \cite{HBK} and used in \cite{CPStepanov}, to a general finite field.
To give a rough description of the method, if $s$ is a positive integer with $16s^2\lessapprox t\lessapprox p/4s$ and $|A+A|\lessapprox st$,
then an auxiliary polynomial $\Phi(X,Y,Z)$ over $\mathbb F_p$ can be constructed so that $\Phi(X,X^t,(1-X)^t)$ is a nonzero polynomial having more zeros (counted with multiplicity) than its degree, a contradiction.
Applying this observation with $s \approx \sqrt{t}$ yields the bound in \eqref{eq:CP}.

In order to guarantee the non-vanishing of $\Phi(X,X^t,(1-X)^t)$ one needs the fact that a nonzero root of a $k$-nomial of degree less than $p$ can have multiplicity at most $k-1$.
For degree $p$ or larger, which we typically encounter over a general finite field, this property fails, for example $X^p+1=(X+1)^p$, a binomial with a zero of multiplicity $p$.
Instead, we exploit this very property, that is, we make heavy use of the Frobenius automorphism.
In our construction we turn the table and construct $\Phi$ so that $\Phi(X,X^t,(1-X)^t)=0$ identically.

Another difference between prime fields and general fields is the additive energy $E^+(A)$ for multiplicative groups.
Heath-Brown and Konyagin \cite{HBK} obtained $E^+(A) \ll |A|^{5/2}$ for any subgroup of a prime field with $t<p^{2/3}$ using the Stepanov method described above;
in \cite{cp2011} it was shown that the constant can be taken to be $16/3$.
Subsequent improvements have been made in the exponent for smaller $t$.
Upper bounds on energy immediately imply lower bounds on $|A+A|$ via the (Cauchy-Schwarz) inequality $|A+A| \ge \frac {t^4}{E^+(A)}$,
so it is natural to ask whether energy estimates extend to $\mathbb F_q$, and can thus be used for estimating $|A+A|$.
Unfortunately they do not.
It is not hard to construct groups with $|A| \approx q^{1/2}$ and $E^+(A)>\frac 14 |A|^3$, as the following example shows.
The example illustrates another key obstruction to generalizing from prime fields to general finite fields: the presence of proper subfields.

\begin{example}
Let $p$ be an odd prime, $q=p^2$ and $A$ be of the form $A=\mathbb F_p^\times \cup \mathbb F_p^\times \omega$, where $\omega\not\in \mathbb F_p$ is a generator for $A$; $t=2(p-1)$.
For any $c$ the number $r(c)$ of solutions to $x+y=c$ with $x,y \in A$ is given by $r(c)=t/2-1$ for $c\in A$,
$r(c)=2$ for $c \not\in A\cup \{0\}$ and $r(0)=t$, as one can readily verify; see ~\cite[proof of Theorem 1.7]{Cochrane2026Fq}.
Thus $|A\cap\F_p^\times|=t/2 \approx q^{1/2}$ and
\[
E^+(A)=\sum_{c} r(c)^2= t(t/2-1)^2+(q-t-1)2^2 +t^2= \tfrac 14 t^3+t^2+t>\tfrac 14 (t+1)^3.
\]
On the other hand, it is plain that $A+A= \mathbb F_q$.

Thus, large energy here coexists with the largest possible sumset.
\end{example}

\begin{remark}
In passing from prime fields to $\mathbb F_q$ with the Stepanov method,
the first thing to note is that the multiplicity of the zeros needed in the construction of $\Phi$ generally exceeds the characteristic $p$,
so instead of using ordinary derivatives one uses Hasse derivatives.
This idea goes back to Schmidt's generalization \cite{schmidt1973, schmidt1976} of Stepanov's original work.
Cochrane and Pinner \cite{cp2006} generalized this approach to an arbitrary mixed exponential sum with rational function arguments.
\end{remark}

\begin{remark}
Stepanov-type approaches have also appeared in recent work studying additive decompositions of multiplicative subgroups;
for example in work of Cochrane, Yip, and Yoo~\cite{CYY}, Hanson and Petridis~\cite{HP}, Kalmynin~\cite{Kalmynin}, Rudnev and Tyrrell~\cite{RT}, and Yip~\cite{Yip, YipRestricted}, among others.
Interestingly, similar obstructions to the two listed earlier also occur for additive decomposition problems when considered in general finite fields instead of prime fields,
often requiring additional non-vanishing hypotheses or restrictions on $|A|$ in terms of $p$.
\end{remark}
\subsection{Outline of Proof}
To describe our construction, write
\[
 |A+A|=t\lambda_A+\varepsilon,\qquad
 \varepsilon=\mathbf1_{\{-1\in A\}},
\]
where $\lambda_A$ is the number of nonzero $A$-cosets in $A+A$.
For the sake of contradiction, suppose that $\lambda_A<\frac 18 t^{1/2}$.
For $u\in A\setminus\{-1\}$, let $x=\frac{u}{1+u}, y=\frac1{1+u}$, so that $x+y=1$ and $x^t=y^t= \eta$, say.
The common value $\eta$ records an {\it occupied} sumset coset, and we call $x,y$ normalized points associated with $\eta$.
We define a linear operator $\mathcal T_j$ on $\overline{\mathbb F}_p[X,Y,Z]$ and construct a polynomial $\Phi(X,Y,Z)$ with specified degrees such that for each common value $\eta$ with associated normalized value $x$,
and $0 \le j < 3\lambda_A$, we have $\mathcal T_j\Phi(X,\eta, \eta)= 0$ and $D^{(j)} \phi(x) = 0$, where $D^{(j)}$ is the $j$-th Hasse derivative, and $\phi(X):=\Phi(X,X^t,(1-X)^t)$.
It will follow that $\phi(X)=0$.
Having constructed $\Phi$, we then select an irreducible factor $H$ of $\Phi$, prove that $H$ in fact has coefficients in $\mathbb F_p$,
and then use $H$ together with its partial derivatives to construct a rational function $R$ with the invariance property $R(u)^{p^r}=u$ on a subset $\mathcal Z$ of $A$ of cardinality at least $\frac 34 t$.
To be precise we establish the following key proposition.

\begin{proposition}
\label{prop:reduction}
Let $q=p^f$, and let $A\le\F_q^\times$ be a multiplicative subgroup of order $t<(\frac 83 q)^{2/3}$.
Suppose that
\[
 |A+A|<t^{3/2}/8.
\]
Then, there exist a positive integer $r$, a nonconstant rational function $R\in\F_p(U)$, and a subset $\mathcal Z\subset A$ with the following properties:
\begin{enumerate}
\item $\deg R<t/5$ and $p^r\deg R<q$.
\item As an identity of rational functions,
\[
 R(-1-U)=-1-R(U).
\]
\item The function $R$ is defined on $\mathcal Z$, and
\[
 |\mathcal Z|>3t/4,\qquad R(u)^{p^r}=u\quad(u\in\mathcal Z).
\]
\end{enumerate}
\end{proposition}
\noindent In fact, we will see that the rational function $R$ is a monomial.
The invariance property in part $(3)$ of the proposition then places $\mathcal Z$ in a proper subfield of $\mathbb F_q$ implying that the intersection of $A$ with the proper subfield exceeds half of $A$, a contradiction.

The bulk of the paper is devoted to proving the proposition, but we start by deducing the main theorem from it.

\section{Proof of Theorem \ref{thm:main}}\label{sec:graph}

In this section we deduce the main theorem from Proposition~\ref{prop:reduction}.
The next two lemmas show that the conclusion of the proposition implies that $R$ is in fact a monic monomial $U^{p^j}$ for some $j \ge 0$.

\begin{lemma}\label{lem:graph}
Let $A\le\F_q^\times$ have order $t$, and let $r\ge0$ be an integer.
Suppose $R\in\F_q(U)$ is a nonzero rational function of degree $e<t$, defined on a subset $\mathcal Z\subset A$, with
\[
 R(u)^{p^r}=u\quad(u\in\mathcal Z),\qquad
 2|\mathcal Z|-t>2e.
\]
Then $R(U)=cU^n$ for some $c\ne0$ and $n\in\mathbb Z$.
\end{lemma}
\begin{proof}
Choose a generator $\gamma$ of $A$, and let $a\in\F_q^\times$ be the unique element with $a^{p^r}=\gamma$.
The subsets $\mathcal Z$ and $\gamma^{-1}\mathcal Z$ of $A$ have an intersection of size at least $2|\mathcal Z|-t>2e$.
For every $u$ in this intersection, we have $u, \gamma u \in \mathcal Z$ and so by assumption
\[
 R(\gamma u)^{p^r}=\gamma u=(aR(u))^{p^r}.
\]
Injectivity of the Frobenius mapping gives $R(\gamma u)=aR(u)$.
Writing $R=P/Q$ in coprime form with $e=\max\{\deg P, \deg Q\}$, we see that each $u \in \mathcal Z \cap \gamma^{-1}\mathcal Z$ is a root of the polynomial
\[
 P(\gamma U)Q(U)-aP(U)Q(\gamma U).
\]
Its degree is at most $2e$, but it has more than $2e$ distinct roots.
It is therefore the zero polynomial, giving
\begin{equation}
 R(\gamma U)=aR(U)
 \quad\text{as rational functions}.\label{eq:generator}
\end{equation}
If $R$ has a zero or pole at $b\in\kk^\times$, equation~\eqref{eq:generator} gives zeros or poles at the $t$ distinct points $b,\gamma b,\ldots,\gamma^{t-1}b$.
This is impossible because $e<t$.
Hence all zeros and poles of $R$ lie at $0$ and $\infty$, and $R(U)=cU^n$.
\end{proof}

\begin{lemma}\label{lem:monomial}
Let $c\in\kk^\times$ and $n\in\mathbb Z\setminus\{0\}$.
Suppose $R(U)=cU^n$ satisfies
\begin{equation}
 R(-1-U)=-1-R(U).\label{eq:equiv}
\end{equation}
Then $R(U)=U^{p^j}$ for some $j\ge0$.
\end{lemma}
\begin{proof}
If $n<0$, the two sides of \eqref{eq:equiv} have their poles at $-1$ and $0$, respectively, which is impossible.
Thus $n>0$, and
\[
 c(-1-U)^n=-1-cU^n.
\]
Comparison of leading coefficients gives $(-1)^n=-1$; evaluation at zero then gives $c=1$.
Hence $(1+U)^n=1+U^n$.
Write $n=p^jm$, where $p\nmid m$.
If $m>1$, the coefficient of $U^{p^j}$ in $(1+U)^n=(1+U^{p^j})^m$ is $m\ne0$, whereas that coefficient in $1+U^n$ is zero.
Therefore $m=1$, as required.
\end{proof}

\begin{proof}[Proof of Theorem~\ref{thm:main} from
Proposition~\ref{prop:reduction}]
Suppose the stated bound fails, and take $r,R,\mathcal Z$ from the proposition.
The function $R$ is nonconstant, and for $e=\deg R<t/5$,
\[
 2|\mathcal Z|-t>t/2>2e.
\]
Lemma~\ref{lem:graph} makes $R$ a monomial, and Lemma~\ref{lem:monomial} gives $R(U)=U^{p^j}$ for some $j\ge0$.

Since $r\ge1$ and $p^j=\deg R$,
\[
 1<p^{r+j}=p^r\deg R<q=p^f.
\]
Thus $0<r+j<f$.
Every $u\in\mathcal Z$ satisfies $u^{p^{r+j}}=u$, so all of $\mathcal Z$ lies in the single proper fixed field
\[
 K=\F_{p^{\gcd(f,r+j)}}<\F_q.
\]
The field is proper because $\gcd(f,r+j)\le r+j<f$.
If $A\cap K^\times=A$, then $\F_p(A)\subseteq K$, contrary to field generation.
Thus $A\cap K^\times$ is a proper subgroup of $A$, has index at least two, and has size at most $t/2$.
This contradicts $\mathcal Z\subset A\cap K^\times$ and $|\mathcal Z|>3t/4$.
\end{proof}

\begin{remark}
The final step above, where we deduce that $A \cap K^\times$ is a proper subgroup of $A$ on account of the assumption $\mathbb F_p(A)=\mathbb F_q$,
is the only place in the proof of Theorem \ref{thm:main} where field generation is needed.
Proposition \ref{prop:reduction} does not require this assumption.
\end{remark}
We proceed to the proof of Proposition \ref{prop:reduction}, starting with some preliminary lemmas.

\section{Preliminary Lemmas}
Let $K$ be any field, $F$ be a rational function over $K$ and $X,T$ indeterminate symbols.
Then in the ring of formal power series $K(X)[[T]]$, $F(X+T)$ admits a formal series expansion
\[
F(X+T)= \sum_{j \ge 0} a_j(X) T^j,
\]
for some $a_j(X) \in K(X)$.
For our interest, we take $K= \overline{\mathbb F}_p$.
The coefficient $a_j(X)$ is called the $j$-th Hasse derivative of $F$, denoted $D^{(j)}F(X)$.
Thus the Hasse derivatives are determined by the relation
\[
 F(X+T)=\sum_{j\ge0}D^{(j)}F(X)T^j .
\]
In particular, $D^{(0)}F(X)=F(X)$, $D^{(1)}F(X) = F'(X)$, the ordinary derivative, and for $j<p$, $D^{(j)}F(X)= F^{(j)}(X)/j!$.
For $j \ge p$ this relationship fails.
Indeed, for any polynomial $F(X)$, if $j \ge p$ then $F^{(j)}(X)=0$ identically, whereas $D^{(j)}F(X) \neq 0$ in general.

\begin{example}
Let $p$ be an odd prime and $F(X)=X^{2p^3}$.
Then over $\mathbb F_p$ we have
\[
F(X+T)= (X+T)^{2p^3}=(X^{p^3}+T^{p^3})^2= X^{2p^3}+2X^{p^3}T^{p^3}+ T^{2p^3},
\]
and so $D^{(j)}F= 0$ unless $j=0,p^3$ or $2p^3$ where it equals $F$, $2X^{p^3}$, $1$ respectively.
\end{example}

If $F$ is regular (no pole) at a point $x \in K$, then $F(x+T) =\sum_{j \ge 0} D^{(j)}F(x) T^j$ in the formal power series ring.
Thus the vanishing of the Hasse derivatives of orders $0,\ldots,m-1$ at $x$ is equivalent to $F$ having a zero of multiplicity at least $m$ at $x$,
that is, $T^m|F(x+T)$, or equivalently $(T-x)^m|F(T)$.
This fact holds in every characteristic.

\begin{lemma}\label{fprimezero}
Let $k$ be a perfect field of characteristic $p$ (such as $\mathbb F_q$), and $f(X)$ be a polynomial/rational function over $k$ such that $f'(X)=0$ identically.
Then $f(X)=F(X^p)$ for some polynomial/rational function $F$ over $k$, and $f(X)=G(X)^p$ for some polynomial/rational function $G$ over $k$.
In particular, if $f(X) \in \mathbb F_p(X)$, then we have $f(X)=F(X^p)=F(X)^p$.
\end{lemma}

\begin{proof}
In both the polynomial and rational function cases, $f(X)$ has a formal Laurent expansion $f(X) = \sum_{i \ge I} a_iX^i$ for some integer $I$ and coefficients $a_i \in k$.
Since $f'(X)=0$ we get $ia_i=0$ for $i \ge I$.
Thus if $p \nmid i$ then $a_i=0$, implying that $f(X)=\sum_{j \ge J} a_{jp} X^{jp}$ for some $J$.
Define $F(X)=\sum_{j \ge J} a_{jp} X^{j}$, so that $f(X)=F(X^p)$.
If $f(X)$ is a polynomial, then plainly $F(X)$ is a polynomial.
Taking $p$-th roots of the coefficients of $F$ in $k$ (which exist since $k$ is perfect), we can write $F(X^p)=G(X)^p$ for some polynomial $G$ over $k$.

If $f(X)$ is a rational function, then the coefficients $a_i$ satisfy a linear recurrence over $k$.
Hence so do the coefficients $a_{jp}$, and so $F(X)$ is a rational function over $k$, say $F(X)=P(X)/Q(X)$ with $P(X),Q(X)$ polynomials.
Taking $p$-th roots of the coefficients of $P,Q$, again we can write $F(X^p)=G(X)^p$ for some rational function $G$.

For the case where $f$ has coefficients in $\mathbb F_p$, it follows that $f(X)=F(X^p)=F(X)^p$, since $a^p=a$ for any $a \in \mathbb F_p$.

\end{proof}

\begin{lemma}
\label{fprimezero2}
Let $k$ be a field of characteristic $p$, and $f(X,Y)$ be a polynomial/rational function over $k$ such that $\frac {\partial f}{\partial X}=0$ and $\frac {\partial f}{\partial Y}=0$ identically.
Then $f(X,Y)=F(X^p,Y^p)$ for some polynomial/rational function $F$ over $k$.
\end{lemma}
\noindent The analogous statement holds for any number of variables.
\begin{proof}
The proof is the same as before so let us just do polynomials.
Write $f=\sum_{i=0}^I \sum_{j=0}^J a_{ij}X^i Y^j$.
The vanishing of $\partial f/\partial X$ gives $a_{ij}i=0$ for all $i,j$, while the vanishing of $\partial f/\partial Y$ gives $a_{ij}j=0$ for all $i,j$.
Thus for any $i,j$, if $a_{ij} \neq 0$ then $p|i$ and $p|j$, and so once again we are left with a polynomial all of whose exponents are multiples of $p$.
\end{proof}
\begin{lemma}\label{lem:implicit}
If $P(Y,Z)\ne0$ is a polynomial over $\overline{\mathbb F}_p$, $t$ is a positive integer with $p \nmid t$, and
\[
 P(X^t,(1-X)^t)=0
\]
identically, then $\deg_Y P\ge t$ and $\deg_Z P\ge t$.
\end{lemma}
\begin{proof}
Let $A$ be the group of $t$-th roots of unity in $\overline{\mathbb F}_p$.
For each $\zeta\in A$, substitute $\zeta X$ for $X$.
Since $\zeta^t=1$, we obtain
\[
 P\bigl(X^t,(1-\zeta X)^t\bigr)=0.
\]
The polynomial $P(X^t,Z)\in\kk(X)[Z]$ is nonzero since the substitution $Y\mapsto X^t$ is injective.
It therefore has the $t$ roots
\[
 (1-\zeta X)^t,\qquad \zeta\in A.
\]
These are distinct, because their coefficients of $X$ are $-t\zeta$ and $p\nmid t$.
Thus $\deg_ZP\ge t$.
Replacing $X$ by $1-X$ in the original relation gives $P((1-X)^t,X^t)=0$.
Applying the same argument to $P(Z,Y)$ gives $\deg_YP\ge t$.
\end{proof}

\section{Proof of Proposition \ref{prop:reduction}: Initial setup and definitions}

Suppose that $A$ is a multiplicative subgroup of $\mathbb F_q^\times$ of order $t$.
Write
\begin{align*}
A+A &= \bigcup_{i=1}^{\lambda_A} (1+u_i)A, \quad \text{with zero added if $-1 \in A$},\\
|A+A|& =t\lambda_A+\varepsilon,\qquad
\varepsilon=\mathbf1_{\{-1\in A\}},
\end{align*}
for some $u_i \in A$.
Now, the mapping $xA \to x^{-t}$ is one-to-one on the cosets of $A$, and so the values $ \eta_i:=(1+u_i)^{-t} \in \mathbb F_q$, $1 \le i \le \lambda_A$, record the distinct cosets in $A+A$.
Put $y_i=(1+u_i)^{-1}$, $x_i=1-y_i= \frac {u_i}{1+u_i}$, so that $x_i+y_i=1$, $x_i^t=y_i^t=\eta_i$.
We call the set $\{\eta_1, \dots, \eta_{\lambda_A}\}$ the set of occupied values.

Going forward, we drop the subscripts above, and simply write for $u\in A\setminus\{-1\}$,
\begin{equation}
x=\frac{u}{1+u},\qquad y=\frac1{1+u},\qquad
x+y=1,\qquad x^t=y^t=\eta.                       \label{eq:normalized}
\end{equation}
We call $\eta$ an {\it occupied} value and $x,y$ {\it normalized points} associated with $\eta$.
There are $\lambda_A$ occupied values $\eta$, and $t-\varepsilon$ distinct normalized pairs $(x,y)$.

In the proposition, we assume that
%$t<q^{2/3}$ and that
$|A+A|<t^{3/2}/8$.
It follows that $t\lambda_A <t^{3/2}/8$ and so $t>64\lambda_A^2 \ge 64$.
Thus we have the standard assumption
\begin{equation} \label{standardassump}
t>64 \lambda_A^2 \ge 64, \quad  \lambda_A<\tfrac 18 \sqrt{t},
\end{equation}
that we can call upon as needed.
In the next section, instead of assuming \eqref{standardassump}, we indicate the requisite assumption in the statement of each lemma.

\section{Constructing the auxiliary polynomial \texorpdfstring{$\Phi(X,Y,Z)$}{Phi(X,Y,Z)}}\label{sec:aux}

%\subsection{The auxiliary polynomial}
Let $A$ be a multiplicative subgroup of $\mathbb F_q^\times$ of order $t$.
For $j \ge 0$, we define a linear operator $\mathcal T_j$ on the set of polynomials $\Phi = \Phi(X,Y,Z)$ over $\overline{\mathbb F}_p$ by
\begin{equation}
\mathcal T_j\Phi=
[X(1-X)]^j[T^j]\,
\Phi\left(X+T,\,
Y\left(\frac{X+T}{X}\right)^t,\,
Z\left(\frac{1-X-T}{1-X}\right)^t\right).           \label{eq:translate}
\end{equation}
Here $[T^j]$ extracts the coefficient of $T^j$ in the formal expansion of the polynomial to its right; the factor $[X(1-X)]^j$ clears the denominators.
For example,
\[
 \mathcal T_1\Phi
 =X(1-X)\tfrac {\partial\Phi}{\partial X}+t(1-X)Y \tfrac {\partial \Phi}{\partial Y}-tXZ\tfrac {\partial \Phi}{\partial Z}.
\]
On the curve $Y=X^t$, $Z=(1-X)^t$, we have
\[
\mathcal T_j\Phi(X,X^t,(1-X)^t)= X^j(1-X)^j[T^j] \Phi(X+T,(X+T)^t,(1-(X+T))^t),
\]
and so with $\phi(X):=\Phi(X,X^t,(1-X)^t)$, the definition of the Hasse derivative gives
\begin{equation} \label{TjDj}
\mathcal T_j\Phi(X,X^t,(1-X)^t)= X^j(1-X)^jD^{(j)}\phi(X), \quad j\ge 0.
\end{equation}

The degrees of the output function $\mathcal T_j \Phi$ are controlled as follows:
\[
\deg_X (\mathcal T_j\Phi) \le j+\deg_X \Phi,\quad  \deg_{Y,Z} (\mathcal T_j \Phi) \le \deg_{Y,Z} \Phi.
\]
Indeed, the contribution of a term $X^iY^vZ^w$ in $\Phi$, to the $T^j$ term in the output $\mathcal T_j\Phi$, is a sum over $a+b=j$, of terms
\[
 (-1)^b\binom{i+tv}{a}\binom{tw}{b}
 X^{i+j-a}(1-X)^{j-b}Y^vZ^w,
\]
with $\deg_X=i+j-a+j-b= i+j \le j+\deg_X \Phi$.
The $Y$ and $Z$ exponents remain the same on every term.

Next, we follow the idea of the Stepanov method used for prime fields, but with key differences.
We construct a polynomial $\Phi$ in the triangular degree region
\begin{equation}\label{eq:triangle-region}
 \deg_X\Phi<3\lambda_A,\qquad \degt\Phi<3\lambda_A.
\end{equation}
such that at each occupied value $\eta$ we have the polynomial identities in $X$
\begin{equation} \label{Tjident1}
 \mathcal T_j\Phi(X,\eta,\eta)=0,\qquad 0\le j<3\lambda_A.
\end{equation}
By \eqref{TjDj} it follows that at each normalized point $x$ we have
\begin{equation} \label{Djphizero}
D^{(j)}\phi(x)=0, \qquad 0 \le j<3\lambda_A.
\end{equation}

The number of coefficients is
\[
 3\lambda_A\binom{3\lambda_A+1}{2}
 =\frac{27\lambda_A^3+9\lambda_A^2}{2}.
\]
For each $j$, the degree of $\mathcal T_j\Phi(X,\eta,\eta)$ is at most $\deg_X \Phi +j \le 3\lambda_A-1+j$ , so setting its coefficients to zero imposes at most $3\lambda_A+j$ conditions.
Altogether there are at most
\[
 \lambda_A\sum_{j=0}^{3\lambda_A-1}(3\lambda_A+j)
 =\frac{27\lambda_A^3-3\lambda_A^2}{2}
\]
homogeneous linear conditions over $\kk$.
The number of coefficients exceeds this count by $6\lambda_A^2>0$, so a nonzero $\Phi$ exists.

In the classical method for prime fields, parameters are chosen so that $\phi(X):=\Phi(X,X^t,(1-X)^t)$ is not identically zero.
This allows one to immediately reach the desired conclusion using the fact that the number of zeros is bounded by the degree.
Here, we have to work harder.
In fact, conditions \eqref{eq:triangle-region} and \eqref{Tjident1} guarantee that in our case, $\phi(X)$ is identically zero.

\begin{lemma}
\label{phizerolem}
If $t>6\lambda_A-1$, then $\phi(X)=0$ identically.
\end{lemma}

\begin{proof}

If $\phi\ne0$, then by \eqref{Djphizero} its number of roots with multiplicity is at least $3\lambda_A(t-\varepsilon)$, whereas by \eqref{eq:triangle-region}, we have
\[
 \deg\phi\le  (3\lambda_A-1)(t+1)<3\lambda_A(t-1),
\]
the latter inequality following from $t>6\lambda_A-1$.
Consequently $\phi=0$.
\end{proof}

\section{The irreducible factor \texorpdfstring{$H$}{H} of \texorpdfstring{$\Phi$}{Phi}}
\label{Hsection}
We assume henceforth that $t>6\lambda_A-1$ so that $\phi(X) =0$.
Factor the nonzero $\Phi(X,Y,Z)$ over $\kk$, and choose an irreducible factor $H=H(X,Y,Z)$ satisfying
\[
H(X,X^t,(1-X)^t)=0.
\]
Such a factor exists because $\kk[X]$ is an integral domain.
By Lemma \ref{lem:implicit}, we have $\deg_X H>0$.
Indeed, if $\deg_X H=0$ then $P(Y,Z):=H(0,Y,Z)=H(X,Y,Z)$ is a nonzero polynomial with $\deg_{Y,Z} P \le \deg_{Y,Z}\Phi<3\lambda_A<t$, satisfying $P(X^t,(1-X)^t)=0$, contradicting the lemma.
It follows that $\deg_{Y,Z} H>0$ as well.
Let $d,B$ denote the source degree and target degree of $H$, so
\begin{equation} \label{Hdeg}
d:=\deg_X H \le 3\lambda_A-1, \quad B:=\deg_{Y,Z} H \le 3 \lambda_A-1.
\end{equation}

\begin{lemma}\label{lem:compare}
Suppose that $t>6\lambda_A-1$ and that $H$ is the irreducible factor defined above.
Suppose $K\ne0$ satisfies $K(X,X^t,(1-X)^t)=0$.
If
\begin{equation} 
 \deg_X H\cdot \deg_{Y,Z}K +\deg_XK\cdot \deg_{Y,Z}H<t,\label{comparison2}
\end{equation}
then $H\mid K$ in $\overline{\mathbb F}_p[X,Y,Z]$.
\end{lemma}

\begin{proof}
Put $C=\deg_{Y,Z}K$, $e=\deg_X K$ so that \eqref{comparison2} becomes $dC+eB<t$.
Use a separate elimination variable $T$, and consider
\[
 \mathcal R(Y,Z)=\Res_T\bigl(H(T,Y,Z),K(T,Y,Z)\bigr).
\]
If $\mathcal R\ne0$, the Sylvester determinant gives
\[
 \deg_Y\mathcal R,\ \deg_Z\mathcal R\le dC+eB<t,
\]
since the determinant has $d$ rows with coefficients of $K$ and $e$ rows with coefficients of $H$.
The resultant also has a B\'ezout identity
\[
 \mathcal R(Y,Z)
 =E(T,Y,Z)H(T,Y,Z)+F(T,Y,Z)K(T,Y,Z)
\]
for some polynomials $E,F\in\kk[T,Y,Z]$.
Substituting $(T,Y,Z)=(X,X^t,(1-X)^t)$, both terms on the right vanish.
Thus
\[
 \mathcal R(X^t,(1-X)^t)=0.
\]
Lemma~\ref{lem:implicit} now rules out a nonzero $\mathcal R$ since $\deg_Y \mathcal R<t$.
Hence $\mathcal R=0$.

If $K$ has positive source degree, then since the leading $T$-coefficients of both $H$ and $K$ are nonzero, the zero resultant gives a common factor of $H,K$ over $\kk(Y,Z)$.
Irreducibility of $H$ implies $H|K$ over $\overline{\mathbb F}_p(Y,Z)$, and it also implies that $H$ is primitive (with a positive source degree) as a polynomial over $\overline {\mathbb F}_p[Y,Z]$.
Thus by Gauss's lemma we have in fact $H\mid K$ in $\kk[X,Y,Z]$.
If $\deg_XK=0$, then $\mathcal R=K^d\ne0$, so this case is excluded.
\end{proof}

Next, we establish uniqueness (up to a constant) of the irreducible $H$ in Lemma \ref{lem:compare}, in the triangular region \eqref{eq:triangle-region}.
We note that the uniqueness includes polynomials that are not factors of $\Phi$.

\begin{lemma}\label{uniquelem}
Suppose that $t>18 \lambda_A^2$ and that $K(X,Y,Z)$ is an irreducible polynomial over $\overline {\mathbb F}_p$ with
\[
\deg_X K \le 3\lambda_A-1, \quad \deg_{Y,Z} K \le 3 \lambda_A-1,
\]
such that $K(X,X^t,(1-X)^t)=0$.
Then $K=cH$ for some $c \in \overline{\mathbb F}_p$.
\end{lemma}
\begin{proof}
Since $H$ satisfies the same degree bounds, the comparison bound \eqref{comparison2} satisfies
\[
 (\deg_X H)\,\degt K+(\deg_XK)\deg_{Y,Z} H 
 \le 2(3\lambda_A-1)^2<18\lambda_A^2<t.
\]
Thus, Lemma~\ref{lem:compare} gives $H\mid K$, and irreducibility of $K$ gives $K=cH$ for some $c\in\kk^\times$.
\end{proof}

We use this uniqueness to show that $H$ has coefficients in $\F_p$.
Let $\sigma$ act on polynomials by applying the Frobenius automorphism to the coefficients and fixing the indeterminates:
\[
 \sigma\left(\sum_{i,v,w}c_{ivw}X^iY^vZ^w\right)
 =\sum_{i,v,w}c_{ivw}^{\,p}X^iY^vZ^w.
\]
Then $\sigma$ is an automorphism of the ring $\overline {\mathbb F}_p [X,Y,Z]$.
Thus, since $H$ is irreducible, so is $\sigma(H)$, and it has the same degrees as $H$.
Since $\sigma$ commutes with the substitution mapping $S(X,Y,Z)=(X,X^t, (1-X)^t)$, we have
\[
 \sigma(H)(X,X^t,(1-X)^t)
 =\sigma\bigl(H(X,X^t,(1-X)^t)\bigr)=0.
\]
Normalize one fixed nonzero coefficient of $H$ to one.
Thus under the hypotheses of Lemma \ref{uniquelem}, uniqueness now gives $\sigma(H)=cH$; comparison at that same coefficient gives $1=c$.
Hence every coefficient of $H$ satisfies $a^p=a$, and so $ H\in\F_p[X,Y,Z]. $ We have established the following.

\begin{lemma}
If $t>18 \lambda_A^2$, then the irreducible factor $H$ can be taken to have coefficients in $\mathbb F_p$, and $H$ is absolutely irreducible over $\mathbb F_p$.
\end{lemma}

\subsection{The two symmetries of the factor \texorpdfstring{$H$}{H}}
\label{subsec:symmetry}

Assume $t>18\lambda_A^2$.
Consider the substitution $(X,Y,Z)\mapsto(1-X,Z,Y)$.
It preserves irreducibility and the degree region.
Since $H(X,X^t,(1-X)^t)=0$ identically, we have
\[
 H(1-X,(1-X)^t,X^t)=0.
\]
Uniqueness therefore gives
\begin{equation}
 H(1-X,Z,Y)=\kappa H(X,Y,Z),\qquad
 \kappa\in\F_p^\times.                         \label{eq:swap}
\end{equation}
Next, consider the substitution
\[
 \iota(X,Y,Z):=(1/X,1/Y,sZ/Y), \quad \quad s=(-1)^t.
\]

\begin{lemma}\label{symlem}
If $t>18\lambda_A^2$, $d=\deg_X H$ and $B=\deg_{Y,Z} H$, then for some $\kappa_\iota\in\F_p^\times$, we have
\begin{equation}
 X^dY^B H(1/X,1/Y,sZ/Y)=\kappa_\iota H(X,Y,Z).              \label{eq:inversion}
\end{equation}
\end{lemma}
\begin{proof}
Clearing denominators, we set
\[
 H_\iota(X,Y,Z):=X^dY^B H(1/X,1/Y,sZ/Y).
\]
Thus at $Y=X^t,\ Z=(1-X)^t$, we have
\[
H_\iota (X,X^t,(1-X)^t)= X^{d+tB}H(\tfrac 1X, \tfrac 1{X^t}, (1-\tfrac 1X)^t)=0,
\]
identically.

Now $\deg_X H_\iota \le d$, and $\deg_{Y,Z} H_\iota \le B$ since a typical monomial $Y^vZ^w$ is sent to $s^wY^{B-v-w}Z^w$, of total degree $B-v\le B$.
Thus $H_\iota$ has comparison bound $2dB<18\lambda_A^2<t$.
Lemma~\ref{lem:compare} gives $H\mid H_\iota$.
Since $\deg_X H_\iota \le \deg_X H$ and $\deg_{Y,Z} H_\iota \le \deg_{Y,Z} H$ we have $H_\iota= \kappa_\iota H$ for some constant $\kappa_\iota$, yielding the identity in the lemma.
\end{proof}

\section{Construction of the rational function \texorpdfstring{$R$}{R} in Proposition \ref{prop:reduction}}
\label{sec:compression}
Assume that $t>18\lambda_A^2$ and that $H \in \mathbb F_p[X,Y,Z]$ is the irreducible factor of $\Phi$ constructed in the preceding section, with degrees satisfying \eqref{Hdeg}.
We now proceed to construct the promised rational function $R(U)$ in Proposition \ref{prop:reduction} from the partial derivatives $\partial H/\partial Y$ and $\partial H/\partial Z$ of $H$.

Write
\begin{equation} \label{hYhZdef}
 h_Y(X):=\tfrac {\partial H}{\partial Y}(X,X^t,(1-X)^t),\qquad
 h_Z(X):=\tfrac {\partial H}{\partial Z}(X,X^t,(1-X)^t).
\end{equation}
Applying the chain rule to $H(X,X^t,(1-X)^t)=0$ gives
\begin{equation} \label{chainrule}
0= \tfrac {\partial H}{\partial X}(X,X^t,(1-X)^t) +tX^{t-1}h_Y(X) -t(1-X)^{t-1}h_Z(X).
\end{equation}

\begin{lemma}\label{lem:nonzero}
Suppose that $t>18\lambda_A^2$.
The polynomials $h_Y$ and $h_Z$ defined in \eqref{hYhZdef} are nonzero.
\end{lemma}
\begin{proof}
Let $d=\deg_X H$, $B=\deg_{Y,Z}H$.
Suppose $h_Y=0$.
If $    \frac {\partial H}{\partial Y} \ne0$, then we obtain a nontrivial relation $\frac {\partial H}{\partial Y}(X,X^t,(1-X)^t)=0$, with
\[
 \deg_X \tfrac {\partial H}{\partial Y} \le d,\qquad \degt \tfrac {\partial H}{\partial Y}\le B-1.
\]
The comparison bound \eqref{comparison2} is at most $d(B-1)+dB<2dB \le 2(3\lambda_A-1)^2<t$.
Lemma~\ref{lem:compare} therefore gives $H\mid \frac {\partial H}{\partial Y } $, contrary to the smaller $(Y,Z)$-degree of $ \frac {\partial H}{\partial Y}$.
We conclude that $ \frac {\partial H}{\partial Y} =0$ as a polynomial in all three variables, not merely after specialization.

Differentiating the swap identity~\eqref{eq:swap} with respect to $Z$ gives
\begin{equation}
\label{swappartial}
 \tfrac {\partial H}{\partial Y}(1-X,Z,Y)=\kappa \tfrac {\partial H}{\partial Z}(X,Y,Z).
\end{equation}
Hence $\frac {\partial H}{\partial Z}=0$ as well and so $h_Z=0$.
The chain rule relation \eqref{chainrule} now reduces to
\[
 \tfrac {\partial H}{\partial X}(X,X^t,(1-X)^t)=0.
\]
If $ \frac {\partial H}{\partial X}\ne0$, its source degree is at most $d-1$ and its total target degree at most $B$.
The comparison bound $dB+(d-1)B<2dB<t$ again gives $H\mid  \frac {\partial H}{\partial X}$, contrary to its smaller source degree.
Thus all three partial derivatives vanish identically, and so by Lemma \ref{fprimezero2}, all exponents of $H$ are therefore divisible by $p$,
and over $\kk$ we can then write $H=G^p$ for some nonconstant polynomial $G$, contradicting irreducibility.
Thus $h_Y\ne0$.
From \eqref{swappartial} we have $h_Y(1-X)=\kappa h_Z(X)$, so $h_Z\ne0$ as well.
\end{proof}

Define
\begin{equation}\label{eq:derivative-polynomials}
 L_Y(X)=X^th_Y(X),\qquad L_Z(X)=(1-X)^th_Z(X),
 \qquad S(X)=\frac{L_Y(X)}{L_Z(X)}.
\end{equation}
The polynomials $L_Y$ and $L_Z$ are specializations of $Y \frac {\partial H}{\partial Y}$ and $Z \frac {\partial H}{\partial Z}$ at $Y=X^t$, $Z=(1-X)^t$.
Multiplication by $Y$ or $Z$ restores the target exponent lowered by differentiation.
We have
\begin{equation}\label{eq:polynomial-degree}
 \deg L_Y,\ \deg L_Z\le d+tB
 \le(3\lambda_A-1)(t+1)<3t\lambda_A.
\end{equation}
Indeed a term $X^uY^vZ^w$ in $H(X,Y,Z)$ contributes a scalar multiple of $X^uX^{tv}(1-X)^{tw}$ to $L_Y$ or $L_Z$, of degree at most $d+t(v+w)$.
The final inequality follows from $t>18\lambda_A^2$.

\subsection{The vanishing of \texorpdfstring{$\mathcal T_jH$}{T\textunderscore j H} on the curve \texorpdfstring{$Y=X^t$}{Y=X\textasciicircum t}, \texorpdfstring{$Z=(1-X)^t$}{Z=(1-X)\textasciicircum t}}

For $j \ge 0$, the function $\mathcal T_j H$ defined in \eqref{eq:translate}, satisfies \eqref{TjDj}, $\mathcal T_jH(X,X^t,(1-X)^t) = X^j(1-X)^jD^{(j)}h(X)$, with $h(X):=H(X,X^t,(1-X)^t)$.
Now $h(X)$ by design is identically zero.
Thus $D^{(j)}h(X)=0$ and
\[
\mathcal T_jH (X,X^t,(1-X)^t)=0
\]
as well, that is, the polynomial $\mathcal T_jH$ vanishes on the curve $Y=X^t$, $Z=(1-X)^t$.
Since $\deg_X(\mathcal T_j H) \le j+\deg_X H$ and $\deg_{Y,Z}(\mathcal T_jH)\le\deg_{Y,Z}H$ we have
\[
\deg_X H \cdot \deg_{Y,Z} \mathcal T_jH + \deg_X \mathcal T_jH \cdot \deg_{Y,Z} H \le  dB+(d+j)B =(2d+j)B,
\]
with $d=\deg_X H$, $B=\deg_{Y,Z}H$, and Lemma~\ref{lem:compare} gives
\begin{equation}
 (2d+j)B<t\quad\Longrightarrow\quad
 \mathcal T_jH(X,Y,Z)=Q_j(X)H                        (X,Y,Z) \label{eq:hasse}
\end{equation}
for some polynomial $Q_j\in\F_p[X]$.
A priori, $Q_j$ is a function of $X,Y,Z$ over $ \overline{\mathbb F}_p$.
To see why $Q_j$ is a function of $X$ alone compare $\deg_{Y,Z}$ of both sides in \eqref{eq:hasse}: $\deg_{Y,Z} \mathcal T_j H \le B$,
while $\deg_{Y,Z} Q_j H = \deg_{Y,Z}Q_j +B$, implying $\deg_{Y,Z}Q_j=0$.
It has coefficients in $\F_p$, since by design, $H$ and consequently $T_jH$, have coefficients in $\mathbb F_p$.
If $\mathcal T_j H =0$ we take $Q_j=0$.

For now, we let $J\ge1$ be any integer satisfying
\begin{equation}\label{eq:admissible-J}
 (2d+J)B<t.
\end{equation}
Then \eqref{eq:hasse} holds for every $1\le j\le J$.
The operators $Y\partial/\partial Y$ and $Z\partial/\partial Z$ commute with $\mathcal T_j$.
Indeed, $\mathcal T_j$ preserves every target exponent pair $(v,w)$, while these two operators multiply its coefficient by $v$ and $w$, respectively.
Since $Q_j$ depends only on $X$, applying them to \eqref{eq:hasse} gives
\[
 \mathcal T_j(Y\tfrac {\partial H}{\partial Y})=Q_j(X)Y \tfrac {\partial H}{\partial Y},\qquad
 \mathcal T_j(Z\tfrac {\partial H}{\partial Z})=Q_j(X)Z \tfrac {\partial H}{\partial  Z}.
\]
Specializing to the curve $Y=X^t$, $Z=(1-X)^t$, using the fact that for any function $F(X,Y,Z)$,
\[
X^j(1-X)^j D^{(j)}F (X,X^t, (1-X)^t)=\mathcal T_jF(X,X^t,(1-X)^t) ,
\]
and applying this with $F=Y\frac {\partial H}{\partial Y}$, $Z \frac {\partial H}{\partial Z}$ therefore gives
\begin{equation}\label{eq:derivative-hasse}
 \begin{aligned}
 [X(1-X)]^jD^{(j)}L_Y&=Q_jL_Y,\\
 [X(1-X)]^jD^{(j)}L_Z&=Q_jL_Z
 \end{aligned}
 \qquad(1\le j\le J).
\end{equation}
Thus the quotients $D^{(j)}L_Y/L_Y$ and $D^{(j)} L_Z/L_Z$ are the same for $1 \le j \le J$, and
\begin{equation} \label{DjLZ}
D^{(j)}L_Y(X) =S(X) D^{(j)} L_Z(X), \quad \quad  1 \le j \le J.
\end{equation}

\subsection{Extracting the root of the quotient}

Applying the Hasse product rule to $L_Y=SL_Z$ and using \eqref{DjLZ} gives
\[
 0=D^{(j)}L_Y-SD^{(j)}L_Z
   =\sum_{k=1}^jD^{(k)}S\,D^{(j-k)}L_Z
 \qquad(1\le j\le J).
\]
For $j=1$, this says $D^{(1)}S\,L_Z=0$, so $D^{(1)}S=0$.
Inductively, the terms with $k<j$ vanish, leaving $D^{(j)}S\,L_Z=0$.
Since $L_Z\ne0$,
\begin{equation}\label{eq:quotient-hasse}
 D^{(j)}S=0\qquad(1\le j\le J).
\end{equation}
Let $r$ be the least integer with $p^r>J$.
Then $r\ge1$ and $1,p,\ldots,p^{r-1}\le J$, and \eqref{eq:quotient-hasse} gives
\[
D^{(1)}S=0, \ D^{(p)}S=0, \dots, D^{(p^{r-1})}S=0.
\]

We claim that this implies $S$ is a perfect $p^r$-th power,
\[
S(X)=G(X)^{p^r}, \quad \quad G(X) \in \mathbb F_p (X).
\]
To begin, $D^{(1)}S=0$ means $S'=0$ and so by Lemma \ref{fprimezero}, $S$ is a $p$-th power, say $S=G_1^p$ with $G_1 \in \mathbb F_p(X)$.
Now, by the definition of the Hasse derivatives, $S(X+T)=\sum_{j\ge 0} D^{(j)}S(X)T^j$ and $G_1(X+T)= \sum_{j\ge 0} D^{(j)}G_1(X)T^j$.
Thus
\[
S(X+T)=G_1(X+T)^p= \sum_{j\ge 0} (D^{(j)}G_1(X))^p T^{jp}. 
\]
Matching the coefficients of $T^p$ in the two expressions for $S(X+T)$ we see that for $r \ge 2$, $0=D^{(p)}S(X)= (D^{(1)}G_1(X))^p$.
Thus $D^{(1)}G_1=0,$ implying that $G_1$ is a $p$-th power, $G_1=G_2^p$ for some $G_2$, $S=G_2^{p^2}$.
If $r \ge 3$ we can repeat this argument with $G_2$ in place of $G_1$ and use $(A+B+\cdots)^{p^2}= A^{p^2}+B^{p^2}+\cdots$ to get $D^{(1)}G_2=0$, whence $S=G_3^{p^3}$ for some $G_3$.
Doing this $r$ times completes the claim.

Put $G=G_r$, so that $S(X)=G(X)^{p^r}$ and define a nonzero rational function $R(U)$,
\begin{equation}\label{eq:R}
  R(U):=G\!\left(\frac{U}{1+U}\right) \in \mathbb F_p(U), 
\end{equation}
so that
\begin{equation} \label{SURU}
S\left(\frac U{1+U}\right) = R(U)^{p^r}.
\end{equation}
Write $G=P/Q$ in coprime form.
Its $p^r$-th power is still in coprime form, with numerator and denominator degrees multiplied by $p^r$.
Also, the fractional linear substitution in \eqref{eq:R} preserves rational degree.
Hence
\begin{equation}\label{eq:degree-general}
 p^r\deg R=\deg S\le d+tB,\qquad
 \deg R\le\frac{d+tB}{J+1},
\end{equation}
where the second bound uses $p^r\ge J+1$.
At this stage $R$ is allowed to be constant; its values at the retained points will establish nonconstancy.

\subsection{Counting zeros and choosing the derivative range}

The same Hasse-derivative identities \eqref{eq:derivative-hasse} control zeros of the original numerator $L_Y$, before any cancellation in $S$.
If $x\notin\{0,1\}$ and $L_Y(x)=0$, then \eqref{eq:derivative-hasse} gives
\[
 D^{(j)}L_Y(x)=Q_j(x)L_Y(x)=0
 \qquad(1\le j\le J).
\]
Each such zero therefore has multiplicity at least $J+1$.
As $L_Y\ne0$,
\begin{equation}\label{eq:zero-general}
 \bigl|\{x\in\kk\setminus\{0,1\}:L_Y(x)=0\}\bigr|
 \le\frac{\deg L_Y}{J+1}\le\frac{d+tB}{J+1}.
\end{equation}
This estimate includes common zeros of $L_Y$ and $L_Z$.

We now choose
\begin{equation}\label{eq:J}
 J=15\lambda_A,
\end{equation}
and note that by the $3\lambda_A$ bound on $d$ and $B$, and assumption $t>63 \lambda_A^2$,
\[
 (2d+J)B<(6\lambda_A+15\lambda_A)\,3\lambda_A
 =63\lambda_A^2<t,
\]
as required by \eqref{eq:admissible-J}.
By \eqref{eq:degree-general}, we have
\begin{equation}\label{eq:small-map}
 \deg R\le\frac{d+tB}{J+1}
 <\frac{3t\lambda_A}{15\lambda_A+1}<\frac t5.
\end{equation}
The number of zeros in \eqref{eq:zero-general} that must be removed from $A$ satisfies the same bound.
Section~\ref{completion} will show that the other values that must be removed total less than $t/20$, leaving more than $3/4$ of the group $A$ for the set $\mathcal Z$ we need to define.

We also derive from \eqref{eq:degree-general} that
\begin{equation}\label{eq:legal}
 \lambda_A<\sqrt t/8, \quad t\le \Big(\frac 83 q\Big)^{2/3}\quad  \Longrightarrow  \quad p^r\deg R\le d+tB<3t\lambda_A
 <\frac38t^{3/2}<q.
\end{equation}
We use this estimate to show that the final fixed field is a proper subfield.

\section{Symmetries of \texorpdfstring{$R(U)$}{R(U)}}\label{sec:map}
In this section we show that $R$ satisfies the symmetry required in the proposition, $R(-1-U)=-1-R(U)$.
Recall from Section \ref{Hsection}, $\iota(X,Y,Z):=(1/X,1/Y,(-1)^tZ/Y)$ and
\[
 X^dY^B(H\circ\iota)=\kappa_\iota H.
\]
Applying $Y\partial/\partial Y$ gives
\[
X^dY^B\Big(Y^{-1} \tfrac {\partial H}{\partial Y }\circ \iota +Y^{-1} (-1)^tZ \tfrac {\partial H}{\partial Z} \circ \iota \Big) = \kappa_\iota \Big(BH-Y\tfrac {\partial H}{\partial Y}\Big),
\]
simplifying to
\[
 X^dY^B\bigl((Y\tfrac {\partial H}{\partial Y}+Z\tfrac {\partial H}{\partial Z})\circ\iota\bigr)=\kappa_\iota(BH-Y\tfrac {\partial H}{\partial Y}).
 \]
On the curve $Y=X^t$, $Z=(1-X)^t$, we have $H =0$ and so the above simplifies further to
\[
X^dY^B\bigl((Y\tfrac {\partial H}{\partial Y}+Z\tfrac {\partial H}{\partial Z})\circ\iota\bigr)
   =-\kappa_\iota Y\tfrac {\partial H}{\partial Y}.
\]
Applying $Z \partial/\partial Z$ gives in a similar manner (both on and off the curve)
\begin{align*}
  X^dY^B\bigl((Z\tfrac {\partial H}{\partial Z})\circ\iota\bigr)
   &=\kappa_\iota Z\tfrac {\partial H}{\partial Z}.
\end{align*}

Dividing the two identities gives for any point on the curve,
\[
\left(\frac {Y{\partial H}/{\partial Y}}{Z {\partial H}/{\partial Z}}\right) \circ \iota +1= - \frac {Y {\partial H}/{\partial Y}}{Z  {\partial H}/{\partial Z}}.
\]
Using the definition
\[
S(X):= \frac {L_Y(X)}{L_Z(X)}= \frac {X^t \frac {\partial H}{\partial Y}(X,X^t,(1-X)^t)}{(1-X)^t \frac {\partial H}{\partial Z}(X,X^t,(1-X)^t)},
\]
the preceding identity simplifies to $S(1/X)+1= -S(X)$, which we write as
\begin{equation}
 S(1/X)=-1-S(X).                               \label{eq:S-inversion}
\end{equation}
The denominator is nonzero by Lemma~\ref{lem:nonzero}, so division is valid in the rational function field.
With $U=\frac X{1-X}$, this transforms via \eqref{SURU} to
\[
R(-1-U)^{p^r} =-1-R(U)^{p^r} = (-1-R(U))^{p^r},
\]
and injectivity gives
\begin{equation}
 R(-1-U)=-1-R(U).           \label{eq:R-equivariance}
\end{equation}
This symmetry is essential for establishing Lemma \ref{lem:monomial}.

\section{Completion of the proof of Proposition \ref{prop:reduction}}\label{completion}

Suppose that $64 \lambda_A^2\le t<(\frac 83q)^{2/3}$.
For each occupied value $\eta$ put $P_\eta(T):=H(T,\eta,\eta)$.
At every normalized point $x=\frac u{1+u}$, $u \in A\setminus\{-1\}$, with that target value ($x^t=(1-x)^t=\eta$), we have
\[
 P_\eta(x)=H(x,x^t,(1-x)^t)=0.
\]
If $P_\eta\ne0$, its degree is at most $d$, so there are at most $d$ such source points $x$.
There are only $\lambda_A$ occupied values $\eta$.
Discarding all such $x$ therefore costs at most $d\lambda_A$ normalized points altogether.
We define
\begin{equation} \label{Zdef}
\mathcal Z:=\{u \in A: u \neq -1, \ P_\eta(T)=0, \ L_Y(x) \neq 0,\  \text{with} \ x=\tfrac u{1+u}, \ \eta=x^t \}.
\end{equation}
Note that the condition, $P_\eta(T)=0$, is stronger than the pointwise equality $P_\eta(x)=0$; it allows us to differentiate the relationship.

Suppose that $u \in \mathcal Z$, $x=\frac u{1+u}$, $\eta=x^t$, so $x$ is a normalized point associated with $\eta$.
In particular $P_\eta(T)=0$ identically.
Differentiating gives for any $x$ associated with $\eta$,
\[
 \tfrac {\partial H}{\partial X}(x,\eta,\eta)=0.
\]
Differentiating the identity $H(X,X^t,(1-X)^t)=0$ gives
\[
\tfrac {\partial H}{\partial X} +tX^{t-1}\tfrac {\partial H}{\partial Y} - t(1-X)^{t-1}\tfrac {\partial H}{\partial Z}=0.
\]
Evaluating at the same point $x$ and using the definitions $y=1-x$, $h_Y(x)=\frac {\partial H}{\partial Y}(x,\eta,\eta)$, $h_Z(x)=\frac {\partial H}{\partial Z}(x,\eta,\eta)$, then gives
\[
 0=t x^{t-1}h_Y(x)-t y^{t-1}h_Z(x)
   =t\eta\left(\frac{h_Y(x)}x-\frac{h_Z(x)}y\right).
\]
Here $x,y,\eta$ are nonzero and $p\nmid t$.
Cancelling $t\eta$ and multiplying by $xy$, we obtain
\begin{equation}\label{eq:verticalderiv}
 yh_Y(x)=xh_Z(x).
\end{equation}

Now, by definition of $\mathcal Z$, $L_Y(x) \neq 0$.
Thus $h_Y(x) \neq 0$ since $L_Y(x)=x^th_Y(x)$.
Equation~\eqref{eq:verticalderiv} then gives $h_Z(x)\ne0$, so $L_Z(x)=y^th_Z(x)\ne0$ as well.
Thus $S$ is defined at $x$.
Its Frobenius root $G$ satisfying $G^{p^r} =S$ is also defined there: a pole of $G$ would remain a pole of $S$.
Consequently $R$ is defined at $u$, and by \eqref{eq:R},
\begin{equation}\label{eq:graphpoints}
 R(u)^{p^r}= G(x)^{p^r}=S(x)
 =\frac{x^th_Y(x)}{y^th_Z(x)}
 =\frac xy=u.
\end{equation}

It remains to count the retained values $u$ in our set $\mathcal Z$.
At most $d\lambda_A$ values fail the condition $P_\eta(T)=0$.
For the other condition, $L_Y(x) \neq 0$, we already saw in \eqref{eq:zero-general} that the number of $x$ that fail this condition is at most $\frac {d+tB}{J+1} \le \frac {3t\lambda_A}{15\lambda_A+1}<\frac t5$, with $J:=15\lambda_A$.

Including the possible parameter $-1$, we thus obtain
\begin{equation}\label{eq:retained-general}
 |\mathcal Z|\ge t-\varepsilon-d\lambda_A
                  -t/5.
\end{equation}
Thus since $d \le 3\lambda_A-1$, we obtain for $t>64\lambda_A^2$,
\[
 \varepsilon+d\lambda_A
 \le1+(3\lambda_A-1)\lambda_A
 \le3\lambda_A^2<\frac{3t}{64}<\frac t{20}.
\]
Therefore with $J=15\lambda_A$, we have
\begin{equation}\label{eq:retained}
 |\mathcal Z|>t-\frac t{20}-\frac t5=\frac{3t}{4}.
\end{equation}
Since $t>64$, this set has more than one point.
The identity $R(u)^{p^r}=u$ makes $R$ injective on $\mathcal Z$, and hence nonconstant.
The degree bounds~\eqref{eq:small-map} and~\eqref{eq:legal}, the symmetry in~\eqref{eq:R-equivariance},
and \eqref{eq:graphpoints}--\eqref{eq:retained} complete the proof of Proposition~\ref{prop:reduction}.

\section{The classical dense-range estimate}\label{sec:dense}

We give the standard Fourier argument for~\eqref{eq:classical-dense}, which complements the main theorem when $t\ge q^{2/3}$.
This estimate requires no field-generation hypothesis.

\begin{proposition}\label{prop:dense}
For every multiplicative subgroup $A\le\F_q^\times$ of order $t$,
\[
 |A+A|\ge \frac{q\,t^3}{t^3+q^2}.
\]
In particular $|A+A|\ge q/2$ when $t\ge q^{2/3}$, and $|A+A| \ge \frac {64}{73} q$ when $t \ge (\frac 83 q)^{2/3}$.
\end{proposition}
\begin{proof}
For an additive character $\psi$ of $\F_q$, write $S(\psi)=\sum_{a\in A}\psi(a)$.
We first show that $|S(\psi)|\le\sqrt q$ for every nontrivial $\psi$.
Put $m=(q-1)/t$.
Multiplicative character orthogonality gives
\[
 S(\psi)=\frac1m\sum_{\chi|_A=1}
     \sum_{z\in\F_q^\times}\chi(z)\psi(z),
\]
where the outer sum is over the $m$ multiplicative characters trivial on $A$.
The inner sum is $-1$ for the trivial character.
For every other $\chi$, its squared absolute value is
\[
 \sum_{u\in\F_q^\times}\chi(u)
       \sum_{y\in\F_q^\times}\psi((u-1)y)
 =(q-1)-\sum_{\substack{u\in\F_q^\times\\u\ne1}}\chi(u)=q.
\]
The triangle inequality now gives $|S(\psi)|\le\sqrt q$.

Additive character orthogonality and Parseval give
\[
 E^+(A)=\frac1q\sum_\psi |S(\psi)|^4
 \le\frac{t^4}{q}+\sum_{\psi\ne1}|S(\psi)|^2
 \le\frac{t^4}{q}+qt.
\]
Finally, Cauchy--Schwarz gives $t^4\le |A+A|E^+(A)$ and hence the stated lower bound.
\end{proof}

\section{Generalization to sums of cosets}\label{cosetsec}
We give a brief description of the modifications required to obtain the estimate for $|A+gA|$ in \eqref{eq:coset}.
Put $\Gamma=gA$, $\alpha=g^t$, and let $\epsilon$ be $1$ if $-1\in\Gamma$, and $0$ otherwise.
The nonzero part of $A+gA$ is a union of $m$ multiplicative $A$-cosets, so $|A+gA|=mt+\epsilon$.
For $u\in\Gamma\setminus\{-1\}$, set
\[
 x=\frac{u}{1+u},\qquad \eta=(1+u)^{-t}.
\]
The $t-\epsilon$ normalized points satisfy $(x^t,(1-x)^t)=(\alpha\eta,\eta)$ and occupy exactly $m$ values of $\eta$.

Suppose $|A+gA|<t^{3/2}/8$.
Since $|A+gA|\ge t$, we have $t>64$, $m\ge1$, and $t>64m^2$.
With $\mathcal T_j$ as before, we impose the condition
\[
 \mathcal T_j\Phi(X,\alpha\eta,\eta)=0
 \quad(0\le j<3m)
\]
for every occupied $\eta$, with $\deg_X\Phi<3m$ and $\deg_{Y,Z}\Phi<3m$.
The coefficient count is $(27m^3+9m^2)/2$ and the number of conditions is at most $(27m^3-3m^2)/2$.
As before, a nonzero $\Phi$ exists and vanishes identically on the curve: otherwise its specialization has at least $3m(t-\epsilon)$ zeros counted with multiplicity but degree at most $(3m-1)(t+1)$.

Choose an absolutely irreducible factor $H$ vanishing on the curve and write $d=\deg_X H$, $B=\deg_{Y,Z}H$.
The lemmas in Section \ref{Hsection} hold as before with $1\le d,B\le3m-1$.
Normalize $H$ over $\F_p$.
Define $h_Y(X),h_Z(X),L_Y(X),L_Z(X)$ as in Section \ref{sec:compression}, and apply the argument of that section with $J=15m$; the required condition is $(2d+J)B<63m^2<t$.

With $\Delta=(d+tB)/(J+1)<t/5$, this gives $P=p^r>J$ and $R\in\F_p(U)$ such that
\[
 \frac{L_Y(X)}{L_Z(X)}
   =R\!\left(\frac{X}{1-X}\right)^P,\qquad
 \deg R\le\Delta,\qquad P\deg R\le d+tB<3mt,
\]
and $R(-1-U)=-1-R(U)$.
The same argument bounds the number of zeros of the original numerator $L_Y$ outside $\{0,1\}$ by $\Delta$, including common zeros before cancellation.

For each occupied $\eta$, either $H(X,\alpha\eta,\eta)$ vanishes identically or it has at most $d$ roots.
Discarding the latter fibres and the numerator zeros, as in Section \ref{completion}, leaves a subset $Z\subseteq\Gamma$ with
\[
 |Z|\ge t-\epsilon-md-\Delta>\frac34t,
 \qquad \epsilon+md\le3m^2<\frac{t}{20}.
\]
At a retained point, $\frac {\partial H}{\partial X}(x,\alpha\eta,\eta)=0$.
Differentiating the curve identity gives $L_Y(x)/x=L_Z(x)/(1-x)$, and hence $R(u)^P=u$ for $u\in Z$.
In particular, $R$ is nonconstant.

The proof of Lemma \ref{lem:graph} applies to $\Gamma$: if $a$ generates $A$, then $|Z\cap a^{-1}Z|\ge2|Z|-t>2\deg R$.
It makes $R$ a monomial and Lemma \ref{lem:monomial} gives $R(U)=U^{p^j}$.
Thus $Z$ lies in the intersection of $\Gamma$ with $K=\F_{p^{\gcd(f,r+j)}}$, and
\[
 p^{r+j}=P\deg R<3mt<\frac38t^{3/2}.
\]
A nonempty intersection $\Gamma\cap K^\times$ is a coset of $A\cap K^\times$.
Its size exceeds $t/2$, so $A\subseteq K$.
Field generation therefore gives $q\le p^{r+j}<3t^{3/2}/8$.

Finally, Fourier Cauchy--Schwarz gives $E(A,gA)\le E(A)$, where $E(B,C)$ counts solutions of $b_1+c_1=b_2+c_2$ with $b_i\in B$, $c_i\in C$, and $E(B)=E(B,B)$.
and so again we have $ |A+gA|\ge\frac{qt^3}{t^3+q^2}>\frac{64}{73}q.\qedhere $

\subsection*{Acknowledgments} ChatGPT by OpenAI was used for brainstorming, computations, drafting, and proofreading during the later stages of the work on this manuscript.
It particularly provided valuable assistance finding counterexamples to earlier considered conjectures,
and also helped identify parts of earlier versions of the argument which were more elaborate than necessary.
Finally, Claude by Anthropic was used for proofreading prior to sharing.

%\nocite{*}

\end{document}